# Indefinite Theta Series on Cones

Martin Westerholt-Raum[1]

**Abstract:** We show that indefinite theta series on cones converge and provide an explicit modular completion. Our completion rests on a convolution of the Gaussian with a piecewise constant function supported on the cone. Our main innovation is to formulate this convolution in terms of euclidean geometry as opposed to hyperbolic geometry. This change of perspective allows us to establish essential asymptotic estimates without further difficulty.

An indefinite theta series is the generating series associated with a vector norm on a suitable cone of an indefinite lattice. As the notion is so flexible, such generating series occur frequently in, e.g. number theory [Kit86; Sie51], infinite dimensional representation theory [BO09; KW01; KW14], enumerative geometry [LZ15; OP06], combinatorics [BO06], and string theory [DMZ12; Gan16]. Indefinite theta series are also connected directly to mock theta series. For decades, mathematicians sought to clarify the mystery around them, which emerged from Ramanujan's Death Bed Letter and his Lost Notebook [Ram00; Ram88]: What are the modular properties of mock theta series? Zwegers resolved it in his pioneering thesis [Zwe02] by offering not one but three explanations, two of which are related to indefinite theta series. In the third section of [Zwe02], Zwegers relates forthright indefinite theta series and mock theta series. Section 1 of his thesis suggests to understand mock theta series in terms of Appell-Lerch sums, which are quotients of indefinite theta series on "cubical cones" by a usual theta series. Zwegers's key discovery was a "modular completion" of indefinite theta series for lattices of signature $(d^+, 1)$ together with an explicit formula for its Fourier expansion. His piece is at the foundation of all applications that we mentioned.

Zwegers's assumption that a lattice has signature $(d^+, 1)$ as opposed to the most general signature $(d^+, d^-)$ with $d^+, d^- \geq 0$ is a major obstacle to several applications. For example, in representation theory of affine superalgebras [KW14; KW16a; KW16b] Zwegers's assumption restricts the number of imaginary roots. In string theory, it limits base spaces for which one can study asymptotics of and exact formulas for topological invariants [AMP13; BM13; Man13].

Given a non-degenerate lattice $L$ in a real quadratic space $V$ with quadratic form $q$ of signature $(d^+, d^-)$ and bilinear form $\langle l, l' \rangle = q(l + l') - q(l) - q(l')$, suitable cones $\mathscr{C} \subseteq V$ yield theta series

$$\theta_{L,0}(\mathscr{C}; \tau, z) = \sum_{l \in L \cap (\mathscr{C} - \frac{v}{y})} \exp\big(2\pi i \,(q(l)\tau + \langle z, l \rangle)\big),$$

$$\tau = x + iy \in \mathbb{H} = \{\tau \in \mathbb{C} : \mathrm{Im}(\tau) > 0\},\ z = u + iv \in V(\mathbb{C}) := V \otimes \mathbb{C}.$$

The key questions are: (1) When does $\theta_{L,0}(\mathscr{C})$ converge? (2) If it converges, what are its modular properties? In other words, what is its modular completion? (3) If $\theta_{L,0}(\mathscr{C})$ is not modular, what is the 1-cocycle attached to it? We answer the first two questions in this work.

---

The author was partially supported by Vetenskapsrådet Grant 2015-04139.

Zwegers's thesis treating the case of signatures $(d^+, 1)$ was published in 2002. The existence of a modular completion of indefinite theta series in general has been announced several times by Zagier and Zwegers in talks over the past years. However, details or explicit formulas are not available. In 2012, the author [Wes15b] extended slightly the approach in Zwegers's thesis to general lattices and cubical cones whose faces meet mutually orthogonally. The orthogonality assumption, which is a major obstacle in applications, was lifted for lattice signatures $(d^+, 2)$ and some cubical cones by Alexandrov, Banerjee, Manschot, and Pioline [Ale+16] in summer 2016. Shortly after [Ale+16] appeared as a preprint, Kudla [Kud16] illustrated how the modular completion of more specific indefinite theta series can be obtained from the theta kernels in [KM86; KM87]. Kudla remarked that the conditions imposed on cones in [Ale+16] are stricter than necessary, and he treated additional cases. After the present manuscript appeared as a preprint, Nazaroğlu [Naz16] extended considerations in [Ale+16] to cubical cones in arbitrary lattices, while still imposing rather strong assumptions on their faces.

Kudla has announced a joint paper with J. Funke in which they generalize [Kud16] to the effect that they interpret $\theta_{L,0}(\mathscr{C})$ in terms of intersection numbers of points with cubes in the Grassmannian of maximal negative definite subspaces of $V$. These cubes correspond to cubical cones in our terminology. Funke and Kudla's ideas yield a modular completion of $\theta_{L,0}(\mathscr{C})$ as an integral of the Kudla-Millson theta kernel. If the cube that corresponds to $\mathscr{C}$ is compact, the integral can be evaluated and one obtains the Fourier expansion of the modular completion. In the non-compact case, the modular completion is merely given as an integral. Using our results for cubical cones, one computes its Fourier expansion as a uniform limit of Fourier expansions in the compact case.

We discuss $\theta_{L,0}(\mathscr{C})$ in more detail. Our focus is on polyhedral cones. We call a cone $\mathscr{C} \subseteq V$ polyhedral, if its cross-section is a (possibly non-compact) polyhedron. Polyhedral cones are sufficient to cover all applications of indefinite theta series that are known to the author. Beyond this they can be used to approximate any other cone. For instance, for Fano varieties the nef cone is polyhedral, and the Morrison-Kawamata cone conjecture extends this in a suitable sense to further varieties. On the other end of the spectrum, the classical Appell-Lerch sum can be related to the lattice with Gram matrix $\left(\begin{smallmatrix} 1 & 1 \\ 1 & 0 \end{smallmatrix}\right)$, which has lattice signature $(1, 1)$. In coordinates, the associated cone is $\mathscr{C} = \{l \in \mathbb{R}^2 : (l_1, l_2 \geq 0) \vee (l_1, l_2 \leq 0)\}$. It is "tetrahedral" and "cubical" at the same time (see Section 2.1), since the projectivization $\mathbb{P}(V)$ of $V$ is one-dimensional. Generalized Appell-Lerch sums can be expanded to indefinite theta series on cubical cones [Ale+16].

To answer the first key question on indefinite theta series on cones, we note that if a polyhedral cone $\mathscr{C}$ contains negative vectors, then $\theta_{L,0}(\mathscr{C})$ does not converge. Consider a polyhedral cone that contains only non-negative vectors and whose isotropic parts are rational. We show in Corollary 4.1 on page 16 that $\theta_{L,0}(\mathscr{C})$ converges to a meromorphic function, whose poles are determined by the isotropic parts of $\mathscr{C}$. Only cones with rational isotropic parts have so-far appeared in applications. As for cases of non-rational isotropic parts of $\mathscr{C}$, convergence of $\theta_{L,0}(\mathscr{C})$ is related to Roth's theorem. These cases are not studied in the present paper, while it might be interesting to examine them with an eye to arithmetic properties of the cone $\mathscr{C}$ itself.

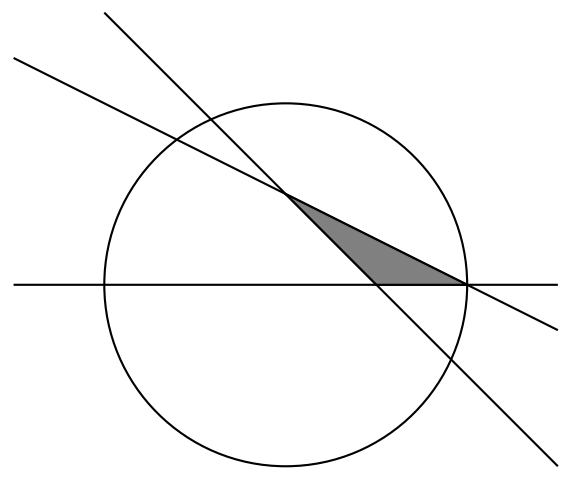

Figure 1: cross section of $\mathscr{C}_{\text{ex}}$

To illustrate our results, we inspect one tetrahedral case as a running example. It is among the simplest cases that are not covered by any prior work.

$$L_{\mathrm{ex}} = \mathbb{Z}^3,\ q_{\mathrm{ex}}(l) = \tfrac{1}{2}(l_1^2 - l_2^2 - l_3^2), \quad \text{and} \qquad \mathscr{C}_{\mathrm{ex}} = \Big\{ l \in \mathbb{R}^3 : \begin{array}{c} (l_3, l_1 - l_2 - 2l_3, 2l_2 + 2l_3 - l_1 \ge 0) \\ \vee\, (l_3, l_1 - l_2 - 2l_3, 2l_2 + 2l_3 - l_1 \le 0) \end{array} \Big\}.$$

Notice that $\mathscr{C}_{\mathrm{ex}}$ contains the isotropic line spanned by $(1\ 1\ 0)$. It is bounded by the orthogonal complements of $(0\ 0\ 1)$, $(1\ 1\ 2)$, and $(1\ 2\ 2)$. Its cross-section for fixed $l_1$ including the circle of isotropic vectors is displayed in Figure 1. The shaded area is the cross section of $\mathscr{C}_{\mathrm{ex}}$ and the displayed lines correspond to its faces.

Recall that the major issue in applications with indefinite theta series is that they are in general not modular. As a remedy Zwegers constructed for lattice signatures $(d^+, 1)$ and polyhedral cones with two faces (i.e. tetrahedral cones) a modular completion of a slight modification of $\theta_{L,0}(\mathscr{C})$. To construct this modification, one first observes that for every polyhedral cone, there is a polynomial $\mathrm{sgn}^+_{\mathscr{C}}(l)$ in sign functions $l \mapsto \mathrm{sgn}\langle l, w\rangle$ with $w \in V$ such that

$$\mathrm{interior}(\mathscr{C}) \subseteq \mathrm{support}(\mathrm{sgn}^+_{\mathscr{C}}) \subseteq \mathrm{closure}(\mathscr{C}) \quad \text{and} \quad \mathrm{sgn}^+_{\mathscr{C}}\big|_{\mathrm{interior}(\mathscr{C})} = \pm 1.$$

For some cones, one can express $\mathrm{sgn}^+_{\mathscr{C}}(l)$ as a linear combination of products of sign functions $\mathrm{sgn}_E(l) := \prod_{w \in E} \mathrm{sgn}\langle l, w\rangle$ with sets $E \subset V$ of size at most $d^-$. For the purpose of this paper, we say that such cones are determined by face indicators. We show in Section 2.1 that tetrahedral and cubical cones are determined by face indicators. In our running example, $\mathrm{sgn}^+_{\mathscr{C}_{\mathrm{ex}}}(l)$ is explicitly given by

$$\tfrac{1}{4}\big(1 + \mathrm{sgn}\langle l, (0\ 0\ 1)\rangle \mathrm{sgn}\langle l, (1\ 1\ 2)\rangle + \mathrm{sgn}\langle l, (0\ 0\ 1)\rangle \mathrm{sgn}\langle l, (1\ 2\ 2)\rangle + \mathrm{sgn}\langle l, (1\ 1\ 2)\rangle \mathrm{sgn}\langle l, (1\ 2\ 2)\rangle\big). \tag{0.1}$$

As a generalization of Zwegers's theta series, we study the meromorphic function

$$\theta_{L,0}(\mathrm{sgn}^+_{\mathscr{C}}; \tau, z) = \sum_{l \in L} \mathrm{sgn}^+_{\mathscr{C}}\big(\sqrt{y}(l + \tfrac{v}{y})\big) \exp\big(2\pi i\, (q(l)\tau + \langle z, l\rangle)\big),$$

Comparing to the definition of $\theta_{L,0}(\mathscr{C})$, the condition that $l \in \mathscr{C} - \frac{v}{y}$ is replaced by $\mathrm{sgn}^+_{\mathscr{C}}(l + \frac{v}{y}) \neq 0$. The difference $\theta_{L,0}(\mathscr{C}) - \theta_{L,0}(\mathrm{sgn}^+_{\mathscr{C}})$ can be expressed in terms of theta series for smaller lattices, if the faces of $\mathscr{C}$ span non-degenerate subspaces of $V$. In that case, the present paper yields a modular completion of $\theta_{L,0}(\mathscr{C})$, which in general has mixed weight.

Zwegers emphasized in many of his talks that the error function that appears in [Zwe02] should be viewed as a smoothened sign. Connecting to this idea, we obtain a smoothened variant $\widehat{\mathrm{sgn}}^+_{\mathscr{C}}$ of $\mathrm{sgn}^+_{\mathscr{C}}$ as a convolution with the Gaussian. The assumption that $\mathscr{C}$ is determined by its face indicators enters. Special cases of $\widehat{\mathrm{sgn}}^+_{\mathscr{C}}$ appear in [Ale+16; Naz16], where they are called generalized error functions. The main novelty in the present manuscript is that our description of $\widehat{\mathrm{sgn}}^+_{\mathscr{C}}$ is based on euclidean geometry as opposed to hyperbolic geometry or contour integration. As a result, we can expose its asymptotic behavior in a straightforward way. In particular, we find that the difference $\widehat{\mathrm{sgn}}^+_{\mathscr{C}} - \mathrm{sgn}^+_{\mathscr{C}}$ decays square-exponentially on every closed cone that does not contain the faces of $\mathscr{C}$.

**Theorem I.** *Let $L$ be an integral lattice. If $\mathscr{C} \subseteq V = L \otimes \mathbb{R}$ is a non-degenerate, non-negative tetrahedral or cubical cone with rational isotropic edges, then*

$$\widehat{\theta}_{L,0}(\widehat{\mathrm{sgn}}^+_{\mathscr{C}}; \tau, z) = \sum_{l \in L} \widehat{\mathrm{sgn}}^+_W\big(\sqrt{y}(l + \tfrac{v}{y})\big) \exp\big(2\pi i\,(q(l)\tau + \langle z, l\rangle)\big)$$

*converges locally to a real-analytic function with meromorphic singularities and transforms like a Jacobi form of suitable level.*

For a discussion of each of the assumptions see the introduction above. Notice that we do not state absolute convergence, which in fact can fail on a certain exceptional set. As a consequence, the Fourier expansion might exhibit "wall-crossing". A more precise statement phrased in terms of vector-valued Jacobi forms is given in Theorem 4.2 on page 16.

*Remark.* Bruinier's and Funke's $\xi$-operator [BF04] plays a decisive role in the modern theory of mock theta series and harmonic weak Maaß forms, which emerged from Zwegers's thesis [Zwe02] and from Bruinier's and Funke's treatment of theta lifts [BF04]. This role can be phrased more systematically in terms of Harish-Chandra modules for real representation theory, as was done in [Sch11; Wes15a]. A description of the analytic properties of $\widehat{\theta}_{L,0}$ along these lines can be inferred from [Wes15b] and the definition Equation (3.1).

To illustrate Theorem I, we return to our running example. Given $v = (\,v_1\ v_2\ v_3\,)$, write

$$\tilde{v}_1 = \tilde{v}_1(v) = v_3 \cdot (\,0\ 0\ 1\,), \quad \tilde{v}_2 = \tilde{v}_2(v) = \tfrac{1}{4}(v_2 + 2v_3 - v_1) \cdot (\,1\ 1\ 2\,), \quad \tilde{v}_3 = \tilde{v}_3(v) = \tfrac{1}{7}(2v_2 + 2v_3 - v_1) \cdot (\,1\ 2\ 2\,).$$

Then $\widehat{\mathrm{sgn}}^+_{\mathscr{C}_{\mathrm{ex}}}(v)$ equals

$$\frac{1}{4} + \frac{1}{4} + \sqrt{3} \int_{\mathbb{R}^2} \frac{\mathrm{sgn}(4v'_1 + 5v'_3)\mathrm{sgn}(5v'_1 + 7v'_3)}{\exp\big(4\pi\,(4(v'_1 - \tilde{v}_1)^2 + 7(v'_3 - \tilde{v}_3)^2 + 10(v'_1 - \tilde{v}_1)(v'_3 - \tilde{v}_3))\big)}\, dv'_1 dv'_3$$
$$+ \sqrt{3} \int_{\mathbb{R}^2} \frac{\mathrm{sgn}(v'_2 + 2v'_3)\mathrm{sgn}(2v'_2 + 7v'_3)}{\exp\big(4\pi\,((v'_2 - \tilde{v}_2)^2 + 7(v'_3 - \tilde{v}_3)^2 + 4(v'_2 - \tilde{v}_2)(v'_3 - \tilde{v}_3))\big)}\, dv'_2 dv'_3.$$

Notice that the second term in the expression (0.1) for $\mathrm{sgn}^+_{\mathscr{C}_{\mathrm{ex}}}$ contributes $\frac{1}{4}$ to $\widehat{\mathrm{sgn}}^+_{\mathscr{C}_{\mathrm{ex}}}(v)$ (defined in (3.9)), as opposed to a general real-analytic function. In this specific case, the reason is that $d^-$ is even and the span of $(\,0\ 0\ 1\,)$ and $(\,1\ 1\ 2\,)$ contains a totally isotropic subspace of dimension $1 = d^- - 1$.

We have intentionally kept this paper as short as possible; We only present material that is strictly necessary to achieve the modular completion of theta series. In Section 1 we recall some basic material. In Section 2 we analyze the geometry of polyhedral cones, and deduce convergence of the incomplete theta series. The main innovation in this paper is that $\widehat{\mathrm{sgn}}^+_{\mathscr{C}}$ is studied in terms of an euclidean geometric decomposition in Section 3. In Section 4 we establish modularity of completed indefinite theta series.

*Acknowledgment* Will be entered after the referees' comments are received.

# 1 Preliminaries

**1.1 Lattices** We let $V$ be a quadratic space over $\mathbb{Q}$ with quadratic form $q$ and scalar product $\langle\cdot,\cdot\rangle$ of signature $(d^+, d^-)$. Its realification is $V(\mathbb{R}) := V \otimes \mathbb{R}$ and its complexification $V(\mathbb{C}) = V \otimes \mathbb{C}$. The set of vectors $v \in V(\mathbb{R})$ with $q(v) > 0$ and $q(v) < 0$ is denoted by $V(\mathbb{R})^+$ and $V(\mathbb{R})^-$, respectively. We write $V(\mathbb{R})^0$ for the set of isotropic vectors in $V(\mathbb{R})$. Given a lattice $L$ in $V$ its dual with respect to the scalar product is denoted by $L^\vee$.

We frequently use the projectivization $\mathbb{P}(V(\mathbb{R})) := (V(\mathbb{R}) \setminus \{0\})/\mathbb{R}$ of $V(\mathbb{R})$. The class in $\mathbb{P}(V(\mathbb{R}))$ that is represented by $v$ is denoted by $[v]$. Given $S \subseteq V(\mathbb{R})$, we write $[S] := \{[v] : v \in S \setminus \{0\}\}$ for the set of projective vectors that are represented by some element in $S$.

The orthogonal complement with respect to $\langle\cdot,\cdot\rangle$ in $V(\mathbb{R})$ of a vector $w$ or a set $W$ is denoted by $w^\perp$ and $W^\perp$, respectively. We write $W \perp W'$ to say that the elements of two set $W$ and $W'$ are mutually orthogonal. Given a definite subspace $U$ and a vector $v$, we write $v_U$ for its orthogonal projection onto $U$. Similarly, if $W$ is a set with definite span, we write $v_W$ for the orthogonal projection onto $\operatorname{span} W$.

**1.2 The Jacobi upper half-space** Let

$$\mathbb{H} = \{\tau \in \mathbb{C} : \operatorname{Im}\tau > 0\}, \qquad \mathbb{H}^{\mathrm{J}}(V) = \mathbb{H} \times V(\mathbb{C}) \tag{1.1}$$

be the Poincaré upper half-plane and the Jacobi upper half-space attached to $V$. We write elements of $\mathbb{H}^{\mathrm{J}}(V)$ as pairs $(\tau, z)$, where the imaginary part of $\tau$ is generally denoted by $y$, and the one of $z$ by $v$. Note that in Section 3, we use the letter $v$ for elements of $V(\mathbb{R})$, but the Jacobi upper half-space will not occur there.

**1.3 The Jacobi groups** We have an action of $\mathrm{SL}_2(\mathbb{R})$ on $\mathbb{H}$ defined by

$$g\tau = \frac{a\tau + b}{c\tau + d}, \quad g = \left(\begin{smallmatrix} a & b \\ c & d \end{smallmatrix}\right) \in \mathrm{SL}_2(\mathbb{R}), \tag{1.2}$$

where the displayed notation for elements of $\mathrm{SL}_2(\mathbb{R})$ is used here and throughout the paper. There is a unique, non-split extension $\mathrm{Mp}_2(\mathbb{R})$ of $\mathrm{SL}_2(\mathbb{R})$ by $\mathbb{Z}/2\mathbb{Z}$, which is called its metaplectic cover. We write $\mathrm{Mp}_2(\mathbb{Z})$ for the preimage of $\mathrm{SL}_2(\mathbb{Z})$ in $\mathrm{Mp}_2(\mathbb{R})$. Elements of $\mathrm{Mp}_2(\mathbb{R})$ can be identified with pairs $(g, \omega)$, where $g \in \mathrm{SL}_2(\mathbb{R})$ and $\omega : \mathbb{H} \to \mathbb{C}$ is a holomorphic square root of $\tau \mapsto c\tau + d$. The action of $\mathrm{SL}_2(\mathbb{R})$ on $\mathbb{H}$ gives rise to one of $\mathrm{Mp}_2(\mathbb{R})$ by the projection $(g, \omega) \mapsto g$.

Fix the Heisenberg group $0 \to \mathbb{R} \to \mathrm{H}(V) \to V(\mathbb{R})^2 \to 0$ with multiplication

$$(\lambda_1, \mu_1, \kappa_1) \cdot (\lambda_2, \mu_2, \kappa_2) = \big(\lambda_1 + \lambda_2,\, \mu_1 + \mu_2,\, \kappa_1 + \kappa_2 + \tfrac{1}{2}\langle\lambda_1, \mu_2\rangle - \tfrac{1}{2}\langle\lambda_2, \mu_1\rangle\big) \tag{1.3}$$

for $\lambda_1, \lambda_2, \mu_1, \mu_2 \in V(\mathbb{R})$ and $\kappa_1, \kappa_2 \in \mathbb{R}$. It carries a right action of $\mathrm{SL}_2(\mathbb{R})$ and $\mathrm{Mp}_2(\mathbb{R})$ by means of the standard representation of $\mathrm{SL}_2(\mathbb{R})$ on the second factor of $V(\mathbb{R})^2 \cong V(\mathbb{R}) \otimes \mathbb{R}^2$.

The real Jacobi group $G^{\mathrm{J}}(V)$ attached to $V$ is the extension

$$0 \longrightarrow \mathrm{H}(V) \longrightarrow G^{\mathrm{J}}(V) \longrightarrow \mathrm{Mp}_2(\mathbb{R}) \longrightarrow 1$$

with multiplication

$$((g_1,\omega_1),h_1)\cdot((g_2,\omega_2),h_2) = \big((g_1,\omega_1)(g_2,\omega_2),\, h_1 g_2 + h_2\big). \tag{1.4}$$

The full Jacobi Group $\Gamma^{\mathrm{J}}(L)$ is the discrete subgroup of $G^{\mathrm{J}}(V)$ generated by $\mathrm{Mp}_2(\mathbb{Z})$ and $L^2 \subset \mathrm{H}(V)$. We usually denoted its elements by $\gamma^{\mathrm{J}}$, identifying them with $(\gamma,\lambda,\mu,\kappa)$ ($\gamma \in \mathrm{Mp}_2(\mathbb{Z})$, $\lambda,\mu \in L$, $\kappa \in \mathbb{Z}$) according to our needs.

**1.4 Weil representations** Assume that $L$ is integral. The Weil representation $\rho_L$ of $L$ is a representation of $\mathrm{Mp}_2(\mathbb{Z})$ with representation space

$$V(\rho_L) := \mathrm{span}_{\mathbb{C}}\{\mathfrak{e}_l : l \in L^\vee/L\}. \tag{1.5}$$

On the generators $T := \big(\left(\begin{smallmatrix}1&1\\0&1\end{smallmatrix}\right), 1\big)$ and $S := \big(\left(\begin{smallmatrix}0&-1\\1&0\end{smallmatrix}\right), \sqrt{\tau}\big)$ of $\mathrm{Mp}_2(\mathbb{Z})$, where $\sqrt{\tau}$ is the principal branch of the holomorphic square root, it is given by

$$\rho_L(T)\,\mathfrak{e}_l := \exp\big(2\pi i\, q(l)\big)\,\mathfrak{e}_l \quad\text{and}\quad \rho_L(S)\,\mathfrak{e}_l := \frac{\sigma(L)}{\sqrt{\# L^\vee/L}} \sum_{l' \in L^\vee/L} \exp\big(-2\pi i\,\langle l, l'\rangle\big)\,\mathfrak{e}_{l'},$$

$$\text{where } \sigma(L) := \frac{1}{\sqrt{\# L^\vee/L}} \sum_{l \in L^\vee/L} \exp\big(-2\pi i\, q(l)\big).$$

The Weil representation extends to a representation of $\Gamma^{\mathrm{J}}(L)$ via the projection $\Gamma^{\mathrm{J}}(L) \twoheadrightarrow \mathrm{Mp}_2(\mathbb{Z})$.

**1.5 Jacobi forms** Assume that $L$ is integral. Let $\phi : \mathbb{H}^{\mathrm{J}}(V) \rightarrow V(\rho_L)$ be a real-analytic function with possible singularities. We call $\phi$ a real-analytic, weak Jacobi form of weight $k$, index $L$, and type $\rho_L$ if for $\gamma^{\mathrm{J}} \in \Gamma^{\mathrm{J}}(L)$ we have

$$\phi\Big(\frac{a\tau+b}{c\tau+d}, \frac{z+\lambda\tau+\mu}{c\tau+d}\Big) = \rho_L(\gamma^{\mathrm{J}})\,\omega(\tau)^{2k} \exp\Big(-2\pi i\Big(\frac{-cq(z+\lambda\tau+\mu)}{c\tau+d} + q(\lambda)\tau + \langle z,\lambda\rangle\Big)\Big)\,\phi(\tau,z).$$

Observe that is suffices to verify that

$$\phi(\tau+1,z) = \rho_L(T)\,\phi(\tau,z), \quad \phi(\tau,z+1) = \phi(\tau,z), \quad\text{and}\quad \phi\Big(\frac{-1}{\tau},\frac{z}{\tau}\Big) = \rho_L(S)\,\sqrt{\tau}^{\,k} \exp\big(2\pi i\, q(z)/\tau\big)\,\phi(\tau,z),$$

in order to check the transformation behavior of $\phi$.

**1.6 Vignéras's modularity result** We write $E = \langle v, \partial_v\rangle$ and $\Delta = \langle \partial_v, \partial_v\rangle$ for the Euler and Laplace differential operators. Note that $\partial_v$ must be viewed as an element of the dual $V^\vee$ of $V$. In other word, in coordinates $v_i$ of $v$, we have $E = \sum_i v_i \partial_{v_i}$ and $\Delta = \sum_{i,j} (m^{-1})_{i,j} \partial_{v_i}\partial_{v_j}$, where $m$ is the Gram matrix associated with $q$.

Vignéras proved a stronger version of the following theorem. It will serve as a major tool for our modularity theorem. We reformulate it slightly, to include the vector valued case, which allows us to work with the full modular group throughout.

**Theorem 1.1 (cf. Vignéras [Vig77]).** *Assume that $q$ is non-degenerate of signature $(d^+, d^-)$. Given a function $h : V(\mathbb{R}) \to \mathbb{C}$ such that $v \mapsto h(v) \exp(-2\pi q(v))$ is a Schwartz function and such that we have*

$$\Delta h / 4\pi = (E - k_0) h \tag{1.6}$$

*for some integer $k_0$, then the theta series*

$$\theta_L(h; \tau, z) := y^{\frac{-k_0}{2}} \sum_{l \in L^\vee} \mathfrak{e}_l \, h\big(\sqrt{y}(l + \tfrac{v}{y})\big) \exp\big(2\pi i (q(l)\tau + \langle z, l \rangle)\big)$$

*converges locally absolutely and uniformly, and is a real-analytic Jacobi form of weight $k_0 + \frac{d^+ + d^-}{2}$, index L, and type $\rho_L$.*

**1.7 One basic estimate** Fix a negative definite subspace $U \subseteq V(\mathbb{R})$, and let $B_D(0) \subset U$ be the ball of radius $D$ with respect to $-q$. To examine the difference of our approximation $\widehat{\mathrm{sgn}}_E$ and the sign-product $\mathrm{sgn}_E$ in Section 3, we need the follow estimate.

**Lemma 1.2.** *Let $U$ be a quadratic space with quadratic form $q$ that is negative definite. If $dv'$ is the measure associated to $-q$, then for every $D \ge 0$, we have*

$$\int_{U \setminus B_D(0)} \exp(4\pi q(v'))\, dv' \Big/ \int_U \exp(4\pi q(v'))\, dv' \ll D^{d-2} \exp(-4\pi D^2).$$

*Proof.* Fix coordinates $v_1, \ldots, v_d$ of $U$. Let $S^{d-1}(r)$ be the $(d-1)$-sphere of radius $r$. Reducing ourselves by a transformation in $\mathrm{GL}(U)$ to the case $q(v) = v_1^2 + \cdots + v_n^2$, we obtain

$$\int_{U \setminus B_D(0)} \exp(4\pi q(v'))\, dv' = \int_D^\infty \exp(-4\pi r^2) \mathrm{vol}(S^{d-1}(r))\, dr = \frac{2\pi^{\frac{d-1}{2}}}{\Gamma(\frac{d-1}{2})} \int_D^\infty r^{d-1} \exp(-4\pi r^2)\, dr$$
$$= \frac{2^{1-d}}{\sqrt{\pi}\Gamma(\frac{d-1}{2})} \int_{4\pi D^2}^\infty r^{\frac{d}{2}-1} \exp(-r)\, dr = \frac{2^{1-d}}{\sqrt{\pi}\Gamma(\frac{d-1}{2})} \Gamma(\tfrac{d}{2}, 4\pi D^2) \ll D^{d-2} \exp(-4\pi D^2). \qquad \blacksquare$$

## 2 Convergence of theta series on cones

Given $C \subseteq \mathbb{P}(V(\mathbb{R}))$, we call

$$\mathscr{C}(C) := \big\{ v \in V(\mathbb{R}) : [v] \in C \big\}$$

the cone associated to $C$. We say that $\mathscr{C}(C)$ is positive, if $C \subset [V(\mathbb{R})^+]$. If $C \subseteq [V(\mathbb{R}) \setminus V(\mathbb{R})^-]$, we say that it is non-negative.

The following theorem is the most elementary convergence theorem for theta series in this paper. We state it separately in order to recurse to it in the proof of the more subtle Theorem 2.7.

**Theorem 2.1.** *Given a closed set $C \subseteq [V(\mathbb{R})^+]$, the theta series $\theta_L(\mathscr{C}(C))$ converges locally absolutely and uniformly.*

*Proof.* Let $\|\cdot\|$ be an auxiliary euclidean norm on $V(\mathbb{R})$. Since $C$ is closed and $\mathbb{P}(V(\mathbb{R}))$ is compact there is a minimum

$$m = \min\big\{q(v) \,:\, v \in V(\mathbb{R}),\, \|v\| = 1,\, [v] \in C\big\},$$

which is positive. We set $\widetilde{q}(v) = m\|v\|^2$, which is a positive definite minorant to $q$ on $\mathscr{C}(C)$. That is, we have $q(v) > m\|v\|^2$ for all $v \in \mathscr{C}(C)$. In particular, we have

$$\Big|\theta_L(\mathscr{C}(C);\tau,z)\Big| \le \sum_{l \in L\cap(\mathscr{C}(C)-\frac{v}{y})} \exp(-2\pi y q(l+\tfrac{v}{y})) \le \sum_{l \in (L+\frac{v}{y})\cap\mathscr{C}(C)} \exp(-2\pi y \widetilde{q}(l)) \le \sum_{v \in L+\frac{v}{y}} \exp(-2\pi y \widetilde{q}(l)).$$

The right hand side is a theta series associated with a positive definite lattice and therefore converges.∎

**2.1 Edges of polyhedral cones** If $C \subseteq \mathbb{P}(V(\mathbb{R}))$ is given as

$$C(W) := \big\{[v] \in \mathbb{P}(V(\mathbb{R})) \,:\, \forall\, w \in W : \langle v, w\rangle \ge 0\big\}$$

for a finite set $W \subset V(\mathbb{R})$, we call $\mathscr{C}(W) := \mathscr{C}(C(W))$ the polyhedral cone with walls $W$. A polyhedral cone is said to be non-degenerate if it has inner points. Throughout this paper we call a cone $\mathscr{C}(W)$ with $\#W = 1 + d^-$ a tetrahedral cone and the corresponding $W$ a tetrahedral set of walls. A set of walls $W$ is called cubical, if $W$ allows for a partition $S$ into $d^-$ pairs such that

$$\operatorname{interior}(C(W)) = \big\{[v] \in \mathbb{P}(V(\mathbb{R})) \,:\, \forall (w,w') \in S : \langle v, w\rangle = \langle v, w'\rangle\big\}.$$

The associated cones are called cubical.

**Lemma 2.2.** *Given a set of walls $W$ such that $\mathscr{C}(W)$ is non-degenerate and non-negative. Then for every $E \subseteq W$ with $\#E > d^-$ we have*

$$V(\mathbb{R})^0 \cap \mathscr{C}(W) \cap E^\perp = \emptyset.$$

*Proof.* The case of $d^- = 0$ is trivial. So let us assume that $d^- > 0$. Suppose that $V(\mathbb{R})^0 \cap \mathscr{C}(W) \cap E^\perp$ was not empty and contained a vector $e$. Choose a positive definite subspace $V' \subset V$ of dimension $d^+ - 1$ that is orthogonal to $e$. Replacing $V$ by the orthogonal complement of $V'$, we can and will reduce our considerations to the case that $d^+ = 1$.

Fix a positive vector $v^+ \in V$. The cross section $A \cap V(\mathbb{R})^0$ of the cone $V(\mathbb{R})^0$ with respect to the affine subspace $A = \{v \in V : \langle v, v^+\rangle = 1\}$ generates $V(\mathbb{R})^0$. Since $\mathscr{C}(W)$ is positive, the analogue holds for its cross section $A \cap \mathscr{C}(W)$. In particular, we may assume that $e \in A$.

Since $d^+ = 1$, $d^- \ge 1$, and $\mathscr{C}(W)$ is polyhedral and non-negative, the cross section $A \cap \mathscr{C}(W)$ is a compact Euclidean polyhedron in $A$. There is a tessellation of $A \cap \mathscr{C}(W)$ by tetrahedrons. Replacing $A \cap \mathscr{C}(W)$ with a tetrahedron containing $e$, we can thus focus on the case of tetrahedral cones. In particular, we can work with the stronger assumption that $e \in V(\mathbb{R})^0 \cap W^\perp$.

Since $\mathscr{C}(W)$ has inner points, there is $v \in V(\mathbb{R})$ with $\langle v, w\rangle > 0$ for all $w \in W$. There exists some $\widetilde{v} \not\in e^\perp$, and by adding a sufficiently small multiple of it to $v$, we can achieve that $\langle v, w\rangle > 0$ for all $w \in W$ and $\langle v, e\rangle \neq 0$. Then for all $t \in \mathbb{R}$, we have $v + te \in \mathscr{C}(W)$, since $e \in W^\perp$. On the other hand, we find that

$$q(v+te) = q(v) + t\langle v, e\rangle + t^2 q(e) = q(v) + t\langle v, e\rangle.$$

For some $t$, we have $q(v+te) < 0$ contradicting the assumption that $\mathscr{C}(W)$ is positive. ∎

For a tetrahedral set $W$, we call $\mathscr{E}(W) = \{E \subseteq W : \#E = d^-\}$ the set of its edges. The orthogonal complements of the $E \in \mathscr{E}(W)$ yield the edges, i.e. the lowest dimensional faces, of the cone $\mathscr{C}(W)$. We call

$$\mathscr{E}(W)^0 := \{E \in \mathscr{E}(W) : E^\perp \cap V(\mathbb{R})^0 \neq \{0\}\}$$

the set of isotropic edges of $W$. In the next two Lemmas, we split up $\mathscr{C}(W)$ into a closed positive set and neighborhoods of isotropic edges.

**Lemma 2.3.** *For every polyhedral cone $\mathscr{C}(W)$ that is non-degenerate and non-negative, we have*

$$\mathscr{C}(W) \cap V(\mathbb{R})^0 = \bigcup_{E \in \mathscr{E}(W)^0} E^\perp = \bigcup_{E \in \mathscr{E}(W)^0} \{v \in V(\mathbb{R}) : \forall w \in E : \langle v, w \rangle = 0\}.$$

*Proof.* As in the proof of Lemma 2.2, it suffices to consider the case of tetrahedral cones and $d^+ = 1$. Fix some $v \in \mathscr{C}(W) \cap V(\mathbb{R})^0$ and let $E$ be a maximal subset of $W$ such that $v \in E^\perp$. We have $E \neq W$ by Lemma 2.2. Suppose that $\#E < d^-$. Then $E^\perp$ contains at least one vector $v'$ with $q(v') < 0$. Observe that $\langle v, w \rangle > 0$ for all $w \in W \setminus E$. For $t \in \mathbb{R}$ with $|t|$ sufficiently small, we therefore have $v + tv' \in \mathscr{C}(W)$. On the other hand,

$$q(v + tv') = q(v) + t\langle v, v' \rangle + t^2 q(v') = t\langle v, v' \rangle + t^2 q(v').$$

If $t \neq 0$, $\operatorname{sgn} t \neq \operatorname{sgn}\langle v, v' \rangle$, and $|t|$ is small enough, then this is negative, contradicting positivity of $\mathscr{C}(W)$. Therefore $\#E = d^-$ as claimed. ■

Given $E \in \mathscr{E}(W)$ fix some auxiliary euclidean norm $\|\cdot\|$ on $V(\mathbb{R})$. With this norm and for $\epsilon > 0$, let

$$E(\epsilon) := E(\|\cdot\|, \epsilon) := \{[v] \in \mathbb{P}(V(\mathbb{R})) : \|v\| = 1, \forall w \in E : 0 \le \langle v, w \rangle < \epsilon\}. \tag{2.1}$$

Notice that we suppress the norm from notation, since for any other choice $\|\cdot\|'$ we have $E(\|\cdot\|', \epsilon') \subseteq E(\|\cdot\|, \epsilon)$ if $\epsilon'$ is small enough. Indeed, $V(\mathbb{R})$ is finite dimensional so that all norms are equivalent to each other. We remark that $E(\epsilon)$ can also be understood in terms of cross sections as at the beginning of the proof of Lemma 2.2.

**Lemma 2.4.** *Given a non-degenerate and non-negative polyhedral cone $\mathscr{C}(W)$, for every $\epsilon > 0$,*

$$C(W) \setminus \bigcup_{E \in \mathscr{E}(W)^0} E(\epsilon) \subseteq \mathbb{P}(V(\mathbb{R}))$$

*is closed and positive.*

*Proof.* Consider the set in the statement. It is closed since $E(\epsilon)$ is open in $C(W)$. It contains no isotropic points by Lemma 2.3, so it is positive. ■

**2.2 Convergence on non-negative polyhedral cones** Consider a non-degenerate and non-negative cone $\mathscr{C}(W)$. For any $E \in \mathscr{E}(W)^0$, write $E(\mathbb{R}) := \mathrm{span}_{\mathbb{R}} E$ and $E(\mathbb{Q}) = E(\mathbb{R}) \cap V$. Write $E(\mathbb{Q})^0$ for the maximal totally isotropic subspace of $V(\mathbb{Q}) \cap E(\mathbb{R})$, which is unique as $E(\mathbb{R})$ contains no negative vectors. Recall that the Witt index of a quadratic space is the dimension of a maximal isotropic subspace. We say that $E \in \mathscr{E}(W)^0$ is a rational edge if the Witt index of $E(\mathbb{R})$ and the dimension of $E(\mathbb{Q})^0$ agree. The next lemma says that we can find neighborhoods of the rational isotropic edges of $C(W)$ in term of polyhedral cones.

**Lemma 2.5.** *Fix a polyhedral, non-degenerate, and non-negative cone $\mathscr{C}(W)$ and an edge $E \in \mathscr{E}(W)^0$ that is rational. Let $E(\mathbb{Q})^- \subset E(\mathbb{Q})$ be a negative definite subspace. There is a set $E' \subset (E(\mathbb{Q})^-)^\perp \cap V(\mathbb{Q})$ of $\dim E(\mathbb{Q})^0$ linearly independent vectors such that $\mathrm{span}\, E'$ is negative definite and $C(E \cap E') \supseteq E(\epsilon)$ for some $\epsilon > 0$.*

*Proof.* Without loss of generality, we can assume that $E(\mathbb{Q})^-$ is trivial by projecting to its orthogonal complement. We can further assume that $E(\mathbb{Q})^0 = E(\mathbb{Q})$ by projecting to the orthogonal complement of a maximal negative definite and rational subspace of $E(\mathbb{R})$. By that orthogonal projection $E(\mathbb{Q})^0$ is mapped to a non-negative subspace of $V(\mathbb{Q})$. This suffices as an assumption on $E(\mathbb{Q})^0$ in the following argument.

Set $n := \dim E(\mathbb{Q})^0$, and observe that we have $d^+, d^- \ge n$. Therefore there is a negative definite, $n$-dimensional subspace $W \subset V(\mathbb{Q})$ such that the pairing $\langle \cdot, \cdot \rangle$ between $E(\mathbb{Q})^0$ and $W$ is non-degenerate. Fix a rational basis $v_1, \ldots, v_n$ of $E(\mathbb{Q})^0$ and a dual basis $w_1, \ldots, w_n$ of $W$. Set $e'_i = e'_i(t) = t w_i - v_i$ for $1 \le i \le n$ and $t > 0$. There is a $t > 0$ such that for any $1 \le i \le d$, we have $q(e'_i) = t^2 q(w_i) - t < 0$, and the matrix $\big(\langle e'_i, e'_j \rangle\big)_{i,j}$ is negative definite. Fix such a $t$, and observe that $\langle v_j, e'_i \rangle \ge 0$ for all $1 \le i, j \le n$. That is, $E' := \big\{ t w_i - v_i \,:\, 1 \le i \le n \big\}$ is a valid choice of $E'$. ■

To make full use of Lemma 2.5, we construct a well-behaved covering of $E(\epsilon)$ by suitable cones $C(E'_i)$.

**Lemma 2.6.** *Let $\mathscr{C}(W)$ be a polyhedral non-degenerate and non-negative cone. Fix $E \in \mathscr{E}(W)^0$ that is rational. Then there is a finite collection of finite sets $E'_i \subset V(\mathbb{R})$ such that for all $\epsilon_i > 0$ we have*

$$E(\epsilon) \subset \bigcup E'_i(\epsilon_i) \quad \text{for some } \epsilon > 0$$

*and the $E'_i$ satisfy:*

*(i) We have $\# E'_i = \dim E'_i(\mathbb{Q}) = d^-$.*

*(ii) We have $E'_i(\mathbb{Q})^0 = E(\mathbb{Q})^0$.*

*(iii) The $E'_i(\epsilon_i)$ are non-negative for sufficiently small $\epsilon_i > 0$.*

*(iv) There is $E'^-_i \subseteq E'_i \cap V(\mathbb{Q})$ such that $\mathrm{span}\, E'^-_i$ is negative definite and $E'_i(\mathbb{Q}) = E(\mathbb{Q})^0 + \mathrm{span}\, E'^-_i$.*

*Proof.* We may replace $V(\mathbb{Q})$ by $(E(\mathbb{Q})^0)^\perp / E(\mathbb{Q})^0$, and thus assume that $E(\mathbb{Q})^0 = \{0\}$.

Choose $\epsilon > 0$ small enough so that $E(\epsilon)$ is positive. We consider $E$ as a point in the space of maximal, framed, negative definite subspaces in $V(\mathbb{R})$. Fix a closed neighborhood $U$ of $E(\epsilon)$ in $[V(\mathbb{R})^+]$.

For any $[v] \in [V(\mathbb{R})^+]$, there is one maximal, framed, negative definite subspace $E'$ and $\epsilon' > 0$ with $E'(\epsilon') \subseteq V(\mathbb{R})^+$ such that $v \in E'(\epsilon')$. Rational points are dense in the space of maximal, framed, negative definite subspaces of $V(\mathbb{R})$. In particular, the compact set $U$ can be covered by the countable set of interiors of $C(E')$ for rational $E'$. Finitely many suffice to cover $C(E) \cap U$, and the set of these $E'$ yields the $E'_i$ in the statement. ∎

**Theorem 2.7.** *Fix a polyhedral cone $\mathscr{C}(W)$ that is non-degenerate and non-negative. If every $E \in \mathscr{E}(W)^0$ is rational, then the theta series $\theta_L(\mathscr{C}(W))$ converges locally absolutely and uniformly on*

$$\big\{(\tau, z) \in \mathbb{H}^{\mathrm{J}}(V) : \forall E \in \mathscr{E}(W)^0 : \emptyset = \mathbb{Z} \cap \langle L + \mathrm{Im}(z)/\mathrm{Im}(\tau), L \cap E(\mathbb{Q})^0 \rangle\big\}.$$

*Proof.* Using notation from (2.1), let $C(\epsilon)^+ = C(W) \setminus \bigcup_{E \in \mathscr{E}(W)^0} E(\epsilon)$ for $\epsilon > 0$. To simplify notation, set $\mathscr{E}(\epsilon) = \mathscr{C}(E(\epsilon))$ and $\mathscr{C}(\epsilon)^+ = \mathscr{C}(C(\epsilon)^+)$. If $\epsilon$ is sufficiently small, we can decompose $\theta_L(\mathscr{C}(W))$ as

$$\theta_L(\mathscr{C}(W); \tau, z) = \theta_L(\mathscr{C}(\epsilon)^+; \tau, z) + \sum_{E \in \mathscr{E}(W)^0} \theta_L(\mathscr{E}(\epsilon); \tau, z).$$

Since $\mathscr{C}(\epsilon)^+$ is positive and closed in $V(\mathbb{R})$ the first summand converges absolutely and uniformly on $\mathbb{H}^{\mathrm{J}}(V)$ by Theorem 2.1. It suffices to demonstrate convergence of $\theta_L(\mathscr{E}(\epsilon))$ for each $E$.

We next employ Lemma 2.6 to reduce our consideration to $E$ such that there is a set of rational vectors $E^- \subset E$ with negative definite $\operatorname{span} E^-$ and $E(\mathbb{Q}) = \operatorname{span} E^- + E(\mathbb{Q})^0$. For such $E$ we obtain from Lemma 2.5 an auxiliary set $E'$ with $\operatorname{span} E'$ negative definite, $E' \perp E^-$, and $C(E \cup E') \supseteq E(\epsilon)$ for sufficiently small $\epsilon$. The pairing of $E(\mathbb{Q})^0$ and $E'(\mathbb{Q})$ is non-degenerate. We may therefore replace $L$ by

$$L \cap E'(\mathbb{Q})^\perp \oplus L \cap E(\mathbb{Q})^0,$$

which has finite index in $L$. In this situation, it suffices to show convergence of

$$\begin{aligned} &\exp(2\pi\, q(\mathrm{Im}(v))) \left|\theta_L(\mathscr{C}(E \cup E'); \tau, z)\right| \\ &\qquad \le \sum_{\substack{l_1 \in L \cap E'(\mathbb{Q})^\perp \\ l_1 + \mathrm{Im}(z_1)/\mathrm{Im}(\tau) \in \mathscr{C}(E^-)}} \;\sum_{\substack{l_2 \in L \cap E(\mathbb{Q})^0 \\ l_2 + \mathrm{Im}(z_2)/\mathrm{Im}(\tau) \in \mathscr{C}(E')}} \exp\big(-2\pi y (q(l_1 + \mathrm{Im}(z)) + \langle l_2, l_1 + \mathrm{Im}(z)\rangle)\big). \end{aligned}$$

As in Zwegers's thesis [Zwe02], for fixed $l_1$ the sum over $l_2$ is bounded, since $\langle l_2, l_1 + v\rangle$ is bounded away from 0. Then the sum over $l_1$ converges, since $\mathscr{C}(E)$ is a positive cone in the orthogonal complement of $E'$. ∎

## 3 Sign functions

We let sgn be the sign function, with $\operatorname{sgn}(0) := 0$. Recall from the introduction that

$$\operatorname{sgn}_E(v) := \prod_{w \in E} \operatorname{sgn}\langle v, w\rangle$$

for any finite subset $E$ of $V(\mathbb{R})$. Notice that $\operatorname{sgn}_E$ is symmetric in the elements of $E$.

**3.1 Smoothened sign function** Assuming that $V$ is negative definite, the normalized convolution

$$\widehat{\mathrm{sgn}}_E(v) := \int_{V(\mathbb{R})} \mathrm{sgn}_E(v')\, \exp(4\pi q(v'-v))\, dv' \Big/ \int_{V(\mathbb{R})} \exp(4\pi q(v'))\, dv' \tag{3.1}$$

is our approximation to $\mathrm{sgn}_E$ (cf. Proposition 3.4). Notice that $\widehat{\mathrm{sgn}}_E(v)$ it smooth and depends only the projection of $v$ to $\mathrm{span}\, E \subseteq V(\mathbb{R})$.

*Remark 3.1.* Up to a change of coordinates and the associated renormalization this smoothened sign is the same as the “higher error function” that appears in Equation (3.44) of [Ale+16] and Equation (4) of [Naz16].

If $V$ is negative semi-definite, then $V(\mathbb{R})^0 \subseteq V(\mathbb{R})$ is a subspace and $V(\mathbb{R})/V(\mathbb{R})^0$ is a negative definite quadratic space. In this situation, we define

$$\widehat{\mathrm{sgn}}_E(v) \;:=\; \widehat{\mathrm{sgn}}_{E^-}(\pi(v))\, \mathrm{sgn}_{E^0}(v), \quad \text{where} \quad E^0 = \{w \in E : q(w) = 0\}, \quad E^- = \{w \in E : q(w) \ne 0\},$$

and $\pi : V(\mathbb{R}) \to V(\mathbb{R})/V(\mathbb{R})^0$ is the canonical projection. For general $V$ and negative semi-definite $\mathrm{span}\, E$, we define

$$\widehat{\mathrm{sgn}}_E(v) := \widehat{\mathrm{sgn}}_E(v_E), \tag{3.2}$$

where the right hand side is computed as with respect to the quadratic space $\mathrm{span}\, E$.

**3.2 Smoothened volume function** To give a decomposition of $\widehat{\mathrm{sgn}}_E$ that exhibits its asymptotic behavior in a precise way, we employ a smoothened version of a volume function. Assuming that $V$ is negative definite, for any finite, non-empty set $E \subset V(\mathbb{R})$, we let

$$S_E(v) := \big\{ v' \in V(\mathbb{R}) : \forall w \in E : \mathrm{sgn}\langle v', w\rangle \ne \mathrm{sgn}\langle v, w\rangle \big\} \tag{3.3}$$

be the sector in $V(\mathbb{R})$ that is opposite to $v$ with respect to the walls corresponding to $E$. Denote its volume by

$$\mathrm{vol}_E(v) := \mathrm{vol}\big(S_E(v)\big), \tag{3.4}$$

where the volume is taken with respect to the quadratic form $-q$ on $V(\mathbb{R})$. Our smoothening

$$\widehat{\mathrm{vol}}_E(v) := \int_{S_E(v)} \exp(4\pi q(v'-v))\, dv' \Big/ \int_{V(\mathbb{R})} \exp(4\pi q(v'))\, dv'$$

depends only on the geometry with respect to $\mathrm{span}\, E$. If $E = \emptyset$, we let $\mathrm{vol}_E(v) = \widehat{\mathrm{vol}}_E(v) = 1$.

**Proposition 3.2.** *For every $v \in V(\mathbb{R})$, we have* $\mathrm{vol}_E(v) = \mathrm{vol}_E(v_E)$.

*Proof.* From its definition, we see that $S_E(v) = S_E(v_E)$. When applying an element of $g \in \mathrm{GL}(E^\perp)$ to $v'$, the normalizing factor $\det(g)$ emerges in both the numerator and the denominator of the defining expression for $\widehat{\mathrm{vol}}_E(v)$, and thus cancels. ∎

The asymptotic properties of $\widehat{\mathrm{vol}}_E$ are described in the next statement.

**Proposition 3.3.** *For $E \neq \emptyset$, we have*

$$\widehat{\mathrm{vol}}_E(v) \ll D^{d-2} \exp\big(-4\pi D^2\big), \qquad D := \mathrm{dist}(v, E^\perp).$$

*Proof.* We observe that the ball $B_D(v)$ does not intersect $S_E(v)$. Therefore, it suffices to estimate

$$\int_{V(\mathbb{R}) \setminus B_D(v)} \exp(4\pi q(v'-v))\, dv' \Big/ \int_{V(\mathbb{R})} \exp(4\pi q(v'))\, dv' = \int_{V(\mathbb{R}) \setminus B_D(0)} \exp(4\pi q(v'))\, dv' \Big/ \int_{V(\mathbb{R})} \exp(4\pi q(v'))\, dv'.$$

This is achieved in Lemma 1.2. ■

**3.3 Decomposition of $\widehat{\mathrm{sgn}}_E$** We write the smoothened sign function as a linear composition of volume functions.

**Proposition 3.4.** *If $V$ is negative definite, we have*

$$\widehat{\mathrm{sgn}}_E(v) = \sum_{E' \subseteq E} 2^{\#E'} \mathrm{sgn}_{E \setminus E'}(v)\, \widehat{\mathrm{vol}}_{E'}(v).$$

*Proof.* We refine the sectors $S_{E'}(v)$ to

$$S^E_{E'}(v) := \big\{ v' \in V(\mathbb{R}) : \forall\, w \in E' : \mathrm{sgn}\langle v', w\rangle \neq \mathrm{sgn}\langle v, w\rangle,\ \forall\, w \in E \setminus E' : \mathrm{sgn}\langle v', w\rangle = \mathrm{sgn}\langle v, w\rangle \big\}, \tag{3.5}$$

the difference being that we now enforce the sign of $\langle v', w\rangle$ also for $w \not\in E'$. Set

$$\widehat{\mathrm{vol}}^E_{E'}(v) = \int_{S^E_{E''}(v)} \exp(4\pi q(v'-v))\, dv' \Big/ \int_{V(\mathbb{R})} \exp(4\pi q(v'))\, dv',$$

and notice that

$$\widehat{\mathrm{vol}}_{E'}(v) = \sum_{E' \subseteq E'' \subseteq E} \widehat{\mathrm{vol}}^E_{E''}(v). \tag{3.6}$$

Define

$$a_n(m) := a_{n-1} - \binom{m}{n-1} a_{n-1}(n-1), \quad m \geq n \geq 1; \quad a_0(m) = (-1)^m. \tag{3.7}$$

Notice that $a_m(m) = (-2)^m$, which can be proved by induction. Indeed, we have

$$a_m(m) = (-1)^m - \sum_{n=0}^{m-1} \binom{m}{n} a_n(n) = (-1)^m - \sum_{n=0}^{m-1} (-2)^n = (-1)^m - (-1)^m + (-2)^m.$$

Fix $v$ and assume that $\mathrm{sgn}_E(v) > 0$. We will show by induction that for every $0 \leq n \leq 1 + d^-$, we have

$$\widehat{\mathrm{sgn}}_E(v) = \sum_{\substack{E' \subseteq E \\ \#E' < n}} (-2)^{\#E'} \widehat{\mathrm{vol}}_{E'}(v) \;+\; \sum_{\substack{E' \subseteq E \\ \#E' \geq n}} a_n(\#E') \widehat{\mathrm{vol}}^E_{E'}(v). \tag{3.8}$$

If $n = 0$, then (3.8) is true by the definitions of $S^E_{E'}$ and $a_0(m)$. We assume that it is true for $n$ and show that it holds for $n \rightsquigarrow n+1$. Rewrite the contribution of $E'$ with $\#E' = n$ to the second summand of (3.8):

$$a_n(n)\widehat{\mathrm{vol}}^E_{E'}(v) = a_n(n)\widehat{\mathrm{vol}}_{E'}(v) \;-\; a_n(n) \sum_{E' \subsetneq E'' \subseteq E} \widehat{\mathrm{vol}}^E_{E''}(v).$$

After inserting this into (3.8), it is straightforward to compare the arising contributions to (3.7).

The case of $\mathrm{sgn}_E(v) < 0$ works analogously. To include the case $\mathrm{sgn}_E(v) = 0$, observe that contributions of $\widehat{\mathrm{vol}}_{E'}$ vanish due to symmetry for $E'$ that satisfy $\langle v, w\rangle = 0$ for some $w \in E \setminus E'$. ∎

**3.4 Degeneration of $\widehat{\mathrm{sgn}}_E$ in families** We consider the behavior of $\widehat{\mathrm{sgn}}_{E(t)}$ for families $E(t)$, $t \in [0,\infty)$. We call such families negative if $\#E(t) = \dim \mathrm{span} E(t)$ is constant on $(0,\infty)$ and $\mathrm{span}\, E(t)$ is negative definite for all $t \in (0,\infty)$.

**Proposition 3.5.** *Assume that $E(t) = E \cup \{w(t)\}$ is a negative family. If $q(w(t)) \to 0$ as $t \to 0$ and $v \not\in w(0)^\perp$, then locally around $v$, we have*

$$\widehat{\mathrm{vol}}_{E(t)}(v) \ll \exp\Big(4\pi \frac{\langle w(t), v\rangle^2}{4q(w(t))}\Big) \to 0 \quad \text{as } t \to 0.$$

*Proof.* This follows when estimating $\mathrm{dist}(v, w^\perp)^2$, $w \in E(t)$ by

$$q\Big(\frac{\langle w(t), v\rangle}{2q(w(t))} w(t)\Big) = \frac{\langle w(t), v\rangle^2}{4q(w(t))},$$

and then applying Proposition 3.3. ∎

**Corollary 3.6.** *With assumptions as in Proposition 3.5, we have locally around $v$*

$$|\widehat{\mathrm{sgn}}_{E(t)}(v) - \widehat{\mathrm{sgn}}_{E(0)}(v)| \ll \exp\Big(\pi \frac{\langle w(t), v\rangle^2}{4q(w(t))}\Big) \quad \text{as } t \to 0.$$

**Proposition 3.7.** *Assume that $E(t) = E^- \cup E^0(t)$ is a negative family, where the Witt index of $E(0)(\mathbb{R})$ equals $\#E^0(t)$. Then $\widehat{\mathrm{sgn}}_{E(t)}$ is a (one-side) differentiable family of smooth functions.*

*Proof.* This can be directly seen from the decomposition in Proposition 3.4. ∎

**3.5 Vignéras's differential equation** We verify Vignéras's differential equation.

**Lemma 3.8.** *Assuming that $V(\mathbb{R})$ is negative definite, we have*

$$\big(2\pi E - \Delta\big)\, \widehat{\mathrm{sgn}}_E = 0.$$

*Proof.* We consider the defining integral of $\widehat{\mathrm{sgn}}_E$. The integrand decays exponentially, so we can differentiate it directly. This shows that the differential $\big(2\pi E - \Delta\big)\widehat{\mathrm{sgn}}_E$ is the integral over a Schwartz function. Now, write $V(\mathbb{R})$ as the union of $U_\epsilon$ and $V(\mathbb{R}) \setminus U_\epsilon$ for an $\epsilon$-neighborhood (with respect to the distance coming from $q$) $U_\epsilon$ of $\{v \in V : \exists w \in E : \langle v, w\rangle = 0\}$. The integrand is differentiable on $V \setminus U_\epsilon$ and $\mathrm{sgn}_E$ is locally constant on there. We therefore see that $(2\pi E - \Delta)\widehat{\mathrm{sgn}}_E(v)$ equals

$$\int_{U_\epsilon} \mathrm{sgn}(v')\big(2\pi E - \Delta\big) \exp(4\pi q(v' - v))\, dv' \;+\; \int_{V(\mathbb{R}) \setminus (U_\epsilon - v)} \big(2\pi E - \Delta\big)\mathrm{sgn}(v' + v) \exp(4\pi q(v'))\, dv'.$$

As $\epsilon \to 0$, this tends to 0, proving the statement. ∎

**3.6 Cones determined by their face indicators** Recall from the introduction that we say that a polyhedral cone $\mathscr{C} = \mathscr{C}(W)$ is determined by its face indicators if there is a linear combination $\mathrm{sgn}^+_{\mathscr{C}}$ of $\mathrm{sgn}_E$ for $E \subset V(\mathbb{R})$ with $\#E \le d^-$ such that

$$\mathrm{interior}(\mathscr{C}) \subseteq \mathrm{support}(\mathrm{sgn}^+_{\mathscr{C}}) \subseteq \mathrm{closure}(\mathscr{C}) \quad \text{and} \quad \mathrm{sgn}^+_{\mathscr{C}}|_{\mathrm{interior}(\mathscr{C})} = \pm 1.$$

The second condition is merely a normalization. By Lemma 2.2, faces of polyhedral cones correspond to sets $E$ as above.

If $\mathscr{C}$ is a cone that is determined by face indicators, we set

$$\widehat{\mathrm{sgn}}^+_{\mathscr{C}}(v) = \sum_E a_E \, \widehat{\mathrm{sgn}}_E(v), \quad \text{where} \quad \mathrm{sgn}^+_{\mathscr{C}}(v) = \sum_E a_E \, \mathrm{sgn}_E(v),\ a_E \in \mathbb{C}. \tag{3.9}$$

**Proposition 3.9.** *Cubical cones are determined by their face indicators.*

*Proof.* Let $\mathscr{C} = \mathscr{C}(W)$ be a cubical cone. By definition there is a partition $S$ of $W$ into $d^-$ pairs such that

$$\mathrm{sgn}^+_{\mathscr{C}}(v) = 2^{-\#W/2} \prod_{(w,w') \in S} \big( \mathrm{sgn}\langle v, w\rangle + \mathrm{sgn}\langle v, w'\rangle \big)$$

satisfies the desired properties. ∎

**Proposition 3.10.** *Tetrahedral cones are determined by their face indicators.*

*Proof.* Given a tetrahedral cone $\mathscr{C} = \mathscr{C}(W)$, define

$$\widetilde{\mathrm{sgn}}^+_{\mathscr{C}}(v) := \prod_{\substack{w,w' \in W \\ w \ne w'}} \big( \mathrm{sgn}\langle v, w\rangle + \mathrm{sgn}\langle v, w'\rangle \big) \cdot \begin{cases} \big((\#W-1)\#W\big)^{-1}, & \text{if } \#W \text{ is odd;} \\ \#W^{-2} \prod \mathrm{sgn}\langle v, w\rangle, & \text{if } \#W \text{ is even.} \end{cases} \tag{3.10}$$

Normalization is chosen in such a way that $\widetilde{\mathrm{sgn}}^+_{\mathscr{C}}(v) = \pm 1$ in the interior of $\mathscr{C}(W)$.

Let $R = \mathbb{C}[s_w : w \in W]$ be the polynomial ring in formal variables indexed by $W$ and $I = R(s_w^2 : w \in W)$ be an ideal. There is a unique polynomial $\widetilde{p} \in R$ such that

$$\widetilde{\mathrm{sgn}}^+_{\mathscr{C}}(v) \;=\; \widetilde{p}\big(s_w \to \mathrm{sgn}\langle v, w\rangle\big).$$

Further, there is a unique lift $p \in R$ of $p \pmod I$ to $R$ that is linear in each of of the variables. We set

$$\mathrm{sgn}^+_{\mathscr{C}}(v) \;=\; p\big(s_w \to \mathrm{sgn}\langle v, w\rangle\big).$$

Checking the total degree of $p$, we conclude that $\mathrm{sgn}^+_{\mathscr{C}}$ is a linear combination of functions $\mathrm{sgn}_E$ for $E \subset V(\mathbb{R})$ of size at most $d^-$.

We have to verify the support of $\mathrm{sgn}^+_{\mathscr{C}}$ to finish the proof. Observe that $\mathrm{sgn}^+_{\mathscr{C}}$ coincides with $\widetilde{\mathrm{sgn}}^+_{\mathscr{C}}$ with possible exceptions on $\cup_w w^\perp$. Given $v$, let $E = \{w \in W \,:\, \langle v, w\rangle = 0\}$. Then it suffices to inspect contributions to $\mathrm{sgn}^+_{\mathscr{C}}$ which do not involve $\mathrm{sgn}\langle v, w\rangle$ for $w \in E$. Collecting these terms in the definition of $\widehat{\mathrm{sgn}}^+_{\mathscr{C}}$, we see that

$$\widehat{\mathrm{sgn}}^+_{\mathscr{C}}(v) := \prod_{\substack{w, w' \in W \\ w \in E \vee w' \in E \\ w \neq w'}} \big(\mathrm{sgn}\langle v, w\rangle + \mathrm{sgn}\langle v, w'\rangle\big)$$

$$\prod_{\substack{w, w' \in W \setminus E \\ w \neq w'}} \big(\mathrm{sgn}\langle v, w\rangle + \mathrm{sgn}\langle v, w'\rangle\big) \cdot \begin{cases} \big((\#W - 1)\#W\big)^{-1}, & \text{if } \#W \text{ is odd;} \\ \#W^{-2} \prod \mathrm{sgn}\langle v, w\rangle, & \text{if } \#W \text{ is even.} \end{cases}$$

If the second factor does not vanish, then $v \in \mathscr{C}(W) \cap E^{\perp}$. The value of the corresponding product in $R/I$ remains the same, since $\mathrm{sgn}\langle v, w\rangle^2 = 1$ for $w \not\in E$. This shows that $v \in E^{\perp} \cap \mathscr{C}(W \setminus E) \subseteq \mathscr{C}(W)$. ∎

# 4 Modular theta series

In this section, we drop notation $v$ for elements in $V(\mathbb{R})$, and write $z = u + iv$ instead. For polyhedral cones $\mathscr{C}$ that are determined by face indicators, we set

$$\theta_L(\mathrm{sgn}^+_{\mathscr{C}};\, \tau, z) = \sum_{l \in L^\vee} \mathfrak{e}_l\, \mathrm{sgn}^+_{\mathscr{C}}\big(\sqrt{y}(l + \tfrac{v}{y})\big) \exp\big(2\pi i\,(q(l)\tau + \langle z, l\rangle)\big). \tag{4.1}$$

As a direct consequence of Theorem 2.7, we find:

**Corollary 4.1.** *Let $\mathscr{C} = \mathscr{C}(W)$ be a polyhedral cone that is non-degenerate, non-negative, and determined by face indicators. If every $E \in \mathscr{E}(W)^0$ is rational, then the theta series $\theta_L(\mathrm{sgn}^+_W)$ converges locally absolutely and uniformly on*

$$\big\{(\tau, z) \in \mathbb{H}^{\mathrm{J}}(V) : \forall E \in \mathscr{E}(W)^0 : \emptyset = \mathbb{Z} \cap \langle L^\vee + \tfrac{v}{y}, L^\vee \cap E(\mathbb{Q})^0\rangle\big\}.$$

For cones determined by their face indicators a completed version of $\theta_L(\mathrm{sgn}^+_{\mathscr{C}})$ arises from $\widehat{\mathrm{sgn}}^+_{\mathscr{C}}$:

$$\theta_L(\widehat{\mathrm{sgn}}^+_{\mathscr{C}};\, \tau, z) := \sum_{l \in L^\vee} \mathfrak{e}_l\, \widehat{\mathrm{sgn}}^+_{\mathscr{C}}\big(\sqrt{y}(l + \tfrac{v}{y})\big) \exp\big(2\pi i(q(l)\tau + \langle z, l\rangle)\big). \tag{4.2}$$

**Theorem 4.2.** *Let $\mathscr{C} = \mathscr{C}(W)$ be a non-degenerate, non-negative tetrahedral or cubical cone such that every $E \in \mathscr{E}(W)^0$ is rational. Then $\theta_L(\widehat{\mathrm{sgn}}^+_{\mathscr{C}}\tau, z)$ with the choice of $\mathrm{sgn}^+_{\mathscr{C}}$ as in Propositions 3.9 and 3.10 converges locally absolutely on*

$$\big\{(\tau, z) \in \mathbb{H}^{\mathrm{J}}(V) : \forall E \in \mathscr{E}(W)^0 : \emptyset = \mathbb{Z} \cap \langle L^\vee + \tfrac{v}{y}, L^\vee \cap E(\mathbb{Q})^0\rangle\big\}.$$

*It is a real-analytic Jacobi form with meromorphic singularities of weight $(d^+ + d^-)/2$, index L, and type $\rho_L$.*

*Proof.* We first establish convergence. To see that (4.2) converges locally absolutely, use the decomposition of $\widehat{\mathrm{sgn}}_E$ in Proposition 3.4 and the asymptotic behavior of $\widehat{\mathrm{vol}}_{E'}$ in Proposition 3.3. In the tetrahedral case, fix $E' \subset W$ with $\#E' < \#W$. Such $E'$ index all possible $\widehat{\mathrm{vol}}_{E'}$ terms in $\widehat{\mathrm{sgn}}^{+}_{W,w}$. The contribution to (4.2) is

$$\sum_{l \in L^\vee} \mathfrak{c}_l \sum_{E \in \mathscr{E}(W)} p_{E'}(l + \tfrac{v}{y})\, \mathrm{sgn}^{+}_{E \setminus E'}(l + \tfrac{v}{y})\, \widehat{\mathrm{vol}}_{E'}\big(\sqrt{y}(l + \tfrac{v}{y})\big) \exp\big(2\pi i (q(l)\tau + \langle z, l \rangle)\big),$$

where $p_{E'}(l)$ is a polynomial in $\mathrm{sgn}\langle v, w\rangle$, $w \in E'$ and $\mathrm{sgn}^{+}_{E\setminus E'}$ is the indicator for the tetrahedral cone spanned by $E \setminus E'$. After using the asymptotic behavior of $\widehat{\mathrm{vol}}_{E'}$ to obtain a partial majorant of $q$, we can invoke Corollary 4.1. Convergence in the cubical case follows along the same lines.

As for modularity, in the case of positive $\mathscr{C}(W)$, Vignéras's modularity theorem suffices, when using Lemma 3.8, because $\widehat{\mathrm{sgn}}^{+}_{\mathscr{C}}(v) \exp(-2\pi q(v))$ is a Schwartz function on $V(\mathbb{R})$. In the general case, we express $\theta_L(\widehat{\mathrm{sgn}}^{+}_{\mathscr{C}})$ as a uniform limit of modular theta series: We first treat the case of isotropic vectors in $E$. This step is completely analogous to [Zwe02]. We let $E(t)$ be a linear approximation of $E(0)$. Then the remainder terms in Corollary 3.6 are sufficient to obtain uniform convergence, when decomposing $V(\mathbb{R})$ as in the proof of Proposition 2.7 in [Zwe02]. To also cover the case of $E(\mathbb{Q})^0$ that is not spanned by isotropic vectors in $E$, we can use Proposition 3.7 directly. ■

---

Chalmers tekniska högskola och Göteborgs Universitet, Institutionen för Matematiska vetenskaper, SE-412 96 Göteborg, Sweden

E-mail: martin@raum-brothers.eu

Homepage: http://raum-brothers.eu/martin